\documentclass{amsart}

\newtheorem{theorem}{Theorem}[section]
\newtheorem{lemma}[theorem]{Lemma}

\theoremstyle{definition}

\theoremstyle{remark}

\numberwithin{equation}{section}

\newcommand\IC{{\mathbb{C}}}
\newcommand\IR{{\mathbb{R}}}

\newcommand\cS{{\mathcal S}}

\newcommand\BH{{{\mathcal B}({\mathcal H})}}

\newcommand\bx{{\bf x}}
\newcommand\by{{\bf y}}

\newcommand\tr{{\rm tr}}

\newcommand\diag{{\rm diag}}

\begin{document}
\openup .3 \jot

\title
{Numerical Range, Dilation, and Completely Positive Maps}

\author{Chi-Kwong Li}
\address{Department of Mathematics, College of William and Mary, 
Williamsbug, Virginia 23187, USA}
\email{ckli@math.wm.edu}

\author{Yiu-Tung Poon}
\address{Department of Mathematics,
Iowa State University, Ames 50011, USA}
\email{ytpoon@iatate.edu}


\date{}

\dedicatory{This paper is dedicated to Professor Tsuyoshi Ando.}


\begin{abstract}
 A proof using the theory of completely positive maps
 is given to the fact that if 
$A \in M_2$, or $A \in M_3$ has a reducing eigenvalue,
then every bounded linear operator $B$ with $W(B) \subseteq W(A)$
has a dilation of the form $I \otimes A$.
This gives a unified treatment for the different cases of
the result obtained by researchers using different techniques. 
\end{abstract}
\maketitle

\section{Introduction}

Let ${\mathcal B}({\mathcal H})$ be the set of bounded linear operators acting 
on a Hilbert space ${\mathcal H}$ with inner product $\langle \bx,\by\rangle$. 
If ${\mathcal H}$ has dimension $n$, we  identify ${\mathcal B}({\mathcal H})$ with
$M_n$ and  ${\mathcal H}=\IC^n$ with $\langle \bx,\by\rangle = \by^*\bx$.  

The numerical range of  $A \in {\mathcal B}({\mathcal H})$ is defined and denoted by
$$W(A) = \{ \langle A\bx,\bx\rangle : \bx \in {\mathcal H}, \langle \bx,\bx\rangle  = 1\}.$$
We say that an operator $B \in {\mathcal B}({\mathcal H})$ admits a dilation 
$A\in {\mathcal B}({\mathcal K})$ if there is a partial isometry $X: {\mathcal H} \rightarrow {\mathcal K}$ such that $X^*X = I_{\mathcal H}$
and $X^*AX = B$.  If $B$ admits a dilation $I_{{\mathcal K}_1}\otimes A\in {\mathcal B}({\mathcal K}_1\otimes {\mathcal K})$
for some Hilbert space ${\mathcal K}_1$,  we will simply say that  
$B$ admits a dilation of the form $I\otimes  A$.

There are interesting connections between the numerical range inclusion and 
dilation relation between two operators. 
For example, the following is known, see 
\cite{M,N}.

\begin{theorem} \label{1.1b} 
Let $A \in M_3$ be a normal matrix.
If $B\in {\mathcal B}({\mathcal H})$ satisfies $W(B)\subseteq W(A) $, then $B$ admits a dilation of the form $I\otimes A$.
\end{theorem}

Also, we have the following
\cite{An,Ar,CL1}.

\begin{theorem} \label{1.1a} 
Let $A \in M_2$.
If $B\in {\mathcal B}({\mathcal H})$ satisfies $W(B)\subseteq W(A) $, then $B$ admits 
a dilation of the form $I \otimes A$.
\end{theorem}

The following theorem was proved in \cite{CL2} that  generalizes   
Theorem \ref{1.1b} and \ref{1.1a}.

\begin{theorem}\label{1.2}
Suppose 
$A \in M_3$ has a non-trivial reducing subspace. If $B \in {\mathcal B}({\mathcal H})$ satisfies  
$W(B) \subseteq W(A)$, then $B$ admits a dilation of the form 
$I \otimes A$.
\end{theorem}

Furthermore, it was shown in \cite{CL1} that the conclusion of Theorem \ref{1.2}
would fail if one considers a general matrix $A \in M_3$ or a normal matrix $A \in M_4$. 
The proofs in \cite{An,Ar,CL1,CL2,M,N} used different techniques. 
In this note, we give a unified proof of the above results.
In particular, we will give a proof of the theorems in Section 1
in an equivalent form involving unital positive and completely  
positive  maps.  We first introduce some background. 

An {\it  operator system} 
$\cS$ of ${\mathcal B}({\mathcal H})$ is a self-adjoint subspace of  ${\mathcal B}({\mathcal H})$ which 
contains $I_{{\mathcal H}}$. A linear map $\Phi:\cS\to {\mathcal B}({\mathcal K})$ is {\it unital} if 
$\Phi\left(I_{{\mathcal H}}\right)= I_{{\mathcal K}}$, $\Phi$ is {\it positive} if 
$\Phi(A)$ is positive semi-definite for every positive semi-definite 
$A\in \cS$,
and  $\Phi$ is said to be {\it completely positive} if 
$I_k\otimes \Phi:M_k(\cS) \to M_k({\mathcal B}({\mathcal K}))$ defined
by $(T_{ij}) \mapsto \left(\Phi(T_{ij})\right)$ is positive 
 for every $k\ge 1$. 

Suppose $A\in {\mathcal B}({\mathcal H})$ and  $B\in {\mathcal B}({\mathcal K})$. Let $\cS$ be the 
operator system spanned by 
$\{I_{{\mathcal H}},A,A^*\}$. Define a unital linear map $\Phi: \cS\to {\mathcal B}({\mathcal K})$  by 
  $\Phi(aI + bA +cA^*) = aI + bB + cB^*$. 
By \cite[Lemma 4.1]{CL2}, $\Phi$ is positive if and only if 
$W(B)\subseteq W(A)$. 
On the other hand,   Stinespring's representation theorem \cite{St} (see also the paragraphs after Theorem 4.1 and Theorem 4.6 in \cite{P}) 
shows that  $\Phi$ is  completely positive if and only if $B$ has a 
dilation of the form $I\otimes A$.  Therefore, Theorems \ref{1.1a} 
-- \ref{1.2} can be stated in the following form.

\begin{theorem} \label{2.1}
Suppose $A = A_0$ or $A_0 \oplus [\mu]$ with $A_0 \in M_2$ and  
$B \in {\mathcal B}({\mathcal H})$.  Define a linear map $\Phi: \cS \rightarrow \BH$ by
$$\Phi(aI + bA +cA^*) = aI + bB + cB^* 
\quad \hbox{ for any } a, b, c \in \IC.$$
Then $\Phi$ is positive if and only if  $\Phi$ is completely positive.
\end{theorem}

\section{Proofs}

To prove Theorem \ref{2.1}, we  need several lemmas, 
some of which are well known. 
In our discussion,  we will let $E_{ij}$ be the basic matrix unit of 
appropriate size, and write $P \ge 0$ if a matrix  or operator 
$P$ is  positive semi-definite.

\begin{lemma} \label{lem2} {\rm (\cite[Corollary 6.7]{P})} Let $\cS$ be an operator system.
Then every positive linear map $\Phi: \cS \rightarrow \BH$
is completely positive for every Hilbert space ${\mathcal H}$ if and only if every positive 
linear map $\Psi: \cS \rightarrow M_n$ is completely 
positive for all positive integer $n$.
\end{lemma}

Recall that $f: \IR^m \rightarrow \IR^m$ is an affine map if it
has the form $\bx \mapsto R\bx + \bx_0$ 
for a real matrix $R \in M_m$ and $\bx_0 \in \IR^m$.
The affine map is invertible if $R$ is invertible, and the 
inverse of $f$ has the form $y \mapsto R^{-1}\by - R^{-1}\bx_0$.
One can extend the definition of affine map to an 
$m$-tuple of self-adjoint operators in $\BH$ by 
$$(A_1, \dots, A_m) \mapsto (A_1, \dots, A_m) (r_{ij} I_{{\mathcal H}})
+ (B_0, \dots, B_m)$$
for a real matrix $R = (r_{ij}) \in M_m$ and an 
$m$-tuple $(B_1, \dots, B_m)$ of self-adjoint operators in $\BH$.
We have the following result which can be easily verified.

\begin{lemma}\label{lem3}
Let $\cS$ be an operator system with a basis
$\{I, A_1, \dots, A_m\}$, and $\Phi: \cS \rightarrow \BH$
  a unital linear map defined by 
$\Phi(A_j) = B_j \in \BH$ for $j = 1, \dots, m$.
Suppose $f$ is an invertible affine map such that 
$f(A_1, \dots, A_m) = (\tilde A_1, \dots, \tilde A_m)$ 
and $f(B_1, \dots, B_m) = (\tilde B_1, \dots, \tilde B_m)$.
Then $\Phi$ is positive (respectively, completely positive)
if and only if the unital map $\tilde \Phi$ defined by
$\tilde \Phi(\tilde A_j) = \tilde B_j$ for $j = 1, \dots, m$,
if positive (respectively, completely positive).
\end{lemma}

\begin{lemma}\label{lem4} Let $\cS = {\rm span}\{E_{jj}: 1 \le j \le m\} \subseteq M_m$. A linear map $\Phi: \cS \rightarrow M_n$
is completely positive if and only if
$\Phi(E_{jj}) \ge 0$ for $j = 1, \dots, m$.
As a result, every positive linear map $\Phi: \cS \rightarrow M_n$
is completely positive.
\end{lemma} 

\it Proof. \rm If $\Phi: \cS \rightarrow M_n$ is positive, then 
$\Phi(E_{jj}) \ge 0$ for all $j = 1, \dots, m$.

Suppose  $\Phi(E_{jj}) \ge 0$ for all $j = 1, \dots, m$.
Let $C = (C_{ij}) \in M_k(\cS)$ be positive semi-definte for a positive
integer $k$.
Then $C = C_{11}\otimes E_{11}+\cdots +C_{mm}\otimes E_{mm}\ge 0$,
where  $C_{jj}\ge 0$ for $j = 1, \dots, m$.  Thus,
$$\left(I_k\otimes \Phi\right)(C)
=C_{11} \otimes \Phi(E_{11}) +\cdots +C_{mm}\otimes \Phi(E_{mm})\ge 0\,.$$
Hence, $\Phi$ is completely positive. \qed

\begin{lemma} \label{lem5} Let 
$\cS = {\rm span}(\{ E_{jj}: 1 \le j \le m\} \cup \{E_{12}+E_{21}\})$ in 
$M_m$ with $m \ge 2$. 
Then every positive linear map  $\Phi:\cS\to M_n$
is completely positive.
\end{lemma}

\it Proof. \rm  Suppose $\Phi: \cS \rightarrow M_n$ is a positive map.
If $m = 2$, the result is due to Choi \cite[Theorem 7]{Choi}.  The proof in \cite{Choi} 
relies on a result of Calderon \cite{C}. 
We give a short and direct proof using basic theory of completely positive
maps for completeness as follows.

Suppose $\Phi$ is positive. Then 
$\Phi(E_{11}), \Phi(E_{22}) \ge 0$, and for any real number $d$, 
\begin{equation} \label{dd}
\Phi(E_{11} + d(E_{12}+E_{21}) + d^2 E_{22})
= \Phi(E_{11}) + d\Phi(E_{12}+E_{21}) + d^2 \Phi(E_{22}) \ge 0.
\end{equation}
Let $C = (C_{ij}) \in M_k(\cS)$ be positive semi-definte for a positive
integer $k$. Then 
$$C  = C_{11} \otimes E_{11} + C_{22} \otimes E_{22}
+ C_{12} \otimes E_{12} + C_{21} \otimes E_{21}$$
such that $C_{21} = C_{12} = C_{21}^*$ so that $C_{12} = C_{21}$
is Hermitian, and 
$$Q = 
\begin{pmatrix} C_{11} & C_{12} \cr C_{21} & C_{22} \cr\end{pmatrix}$$
is positive semi-definite.
We need to show that
\begin{equation} \label{cij}
(I_k \otimes \Phi)(C) = 
C_{11} \otimes \Phi(E_{11}) + C_{12}\otimes \Phi(E_{12}+E_{21})+
C_{22}\otimes \Phi(E_{22}) \ge 0.
\end{equation}
We focus on the case when $C_{11}$ is invertible. The general case  can be derived by continuity argument.
We may replace $C_{ij} = C_{11}^{-1/2} C_{ij} C_{11}^{-1/2}$ for $i,j \in \{1,2\}$ in 
(\ref{cij}), and assume that $C_{11} = I$. Because $C_{12} = C_{12}^*$, we may further 
replace $C_{ij}$ by $U^*C_{ij}U$ in (\ref{cij}) for $i,j \in \{1,2\}$, and assume that   
$C_{11} = I$ and 
$C_{12} = C_{21} = D = \diag(d_1, \dots, d_k)$ with
$d_1, \dots, d_k \in \IR$.
Since $Q \ge 0$, we have
$$\tilde C_{22} = C_{22} - C_{21} C_{11}^{-1}C_{12} 
= C_{22} - D^2 \ge 0.$$
As $\phi(E_{22})$ and $\tilde C_{22}$ are positive semi-definite, 
it follows from (\ref{dd}) that 
\begin{eqnarray*}
&&C_{11} \otimes \Phi(E_{11}) + C_{12}\otimes \Phi(E_{12}+E_{21})+
C_{22}\otimes \Phi(E_{22}) \\
&=& I \otimes \Phi(E_{11}) + D\otimes \Phi(E_{12}+E_{21}) + D^2 \otimes \phi(E_{22})
+ \tilde C_{22} \otimes \Phi(E_{22}) \ge 0
\end{eqnarray*}
as asserted.

\medskip

Suppose $m > 2$ and $k$ is a positive integer.
Let $C = (C_{ij}) \in M_k(\cS)$ be positive semi-definte. Then 
$$C  = C_{11} \otimes E_{11} + C_{22} \otimes E_{22}
+ C_{12} \otimes E_{12} + C_{21} \otimes E_{21} +
\sum_{j=3}^m C_{jj}\otimes E_{jj},$$
where $\sum_{1 \le r,s \le 2} (C_{rs} \otimes E_{rs}) \ge 0$ and
$C_{jj}\ge 0$ for all  $3\le j\le m$. 
We need to show that
$$\left(I_k\otimes \Phi\right)(C)
=\sum_{1 \le r,s \le 2}C_{rs} \otimes \Phi(E_{rs})
+ \sum_{j=3}^m C_{jj} \otimes \Phi(E_{jj})$$
is positive semi-definite.
Since
$\sum_{j=3}^m C_{jj} \otimes \Phi(E_{jj}) \ge 0$, it suffices 
to prove that
$$\sum_{1 \le r,s \le 2}C_{rs} \otimes \Phi(E_{rs}) \ge 0,$$
which is true because  the restriction of 
$\Phi$ to $\{E_{12}, E_{21}, E_{12}+E_{21}\}$ is positive,
and is completely positive by the result
when $m = 2$. The asserted result follows.
\qed

\medskip
Now, we are ready to present the following.

\medskip\noindent
\it Proof of Theorem \ref{2.1}.  \rm By Lemma \ref{lem2}, it suffices to show that if 
$\Phi:\cS \rightarrow M_n$ is 
positive then it is completely positive, where 
$\cS = {\rm span}\{I, A, A^*\}$ satisfies the assumption of the theorem.
We assume that $A$ is not a scalar matrix to avoid trivial consideration.

We will use the fact  that the conclusion will not change
if replace $A$ by $\alpha I + \beta U^*AU$ for any
unitary matrix $U$, 
and $\alpha, \beta \in \IC$ with $\beta \ne 0$.
Furthermore, a unital linear map $\Phi: \cS \rightarrow M_n$ is 
 positive  if and only if 
$B = \Phi(A)$ satisfies $W(B) \subseteq W(A)$ \cite[Lemma 4.1]{CL2}.
Applying an affine transform to $A = A_1+iA_2$
with $(A_1, A_2) = (A_1^*, A_2^*)$ will 
always mean  applying a (real) affine transform to $(A_1,A_2)$.

First, suppose $A$ is normal. If $A \in M_2$ has eigenvalues
$a_1, a_2$,  we may assume that $A = \diag(a_1, a_2)$.
By Lemma \ref{lem3}, we may apply an  invertible affine map to $A$
and assume that  $(a_1,a_2) = (1,0)$. 
Then $\cS = {\rm span}\{E_{11}, E_{22}\}$.
By Lemma \ref{lem4}, every unital
positive linear map  $\Phi: \cS \rightarrow M_n$ is completely positive.

Suppose $A \in M_3$ is normal.
By Lemma \ref{lem3}, we can 
apply  an  invertible affine map and 
assume that  (1) 
$A = \diag(1,r,0)$ with $r \in [0,1]$ if the three eigenvalues
of $A$ are collinear, or (2) 
$A = \diag(1,i ,0)$ otherwise.

Suppose (1) holds, and $\Phi: \cS\rightarrow M_n$
is a unital  positive linear map with
$\Phi(A) = B$. Then $\Phi$ is positive 
if and only if $\bx^*B\bx \in W(A) = [0,1]$ 
for all unit vector $\bx \in \IC^n$, i.e., 
$B$ is a positive semi-definite contraction. 
One can extend $\Phi$ to a map  from $ \hat \cS = 
{\rm span}\{E_{11}, E_{22}, E_{33} \}\subseteq M_3$ to $M_n$ with 
$ \Phi(E_{11}) = B,   \Phi(E_{22}) = 0,   \Phi(E_{33}) = I-B$.
By Lemma \ref{lem4}, $\Phi$ is completely positive on $\tilde \cS$,
and therefore so must be the restriction map on $\cS$.

If (2) holds, then $\cS = {\rm span}\{E_{11}, E_{22}, E_{33}\}$. 
By Lemma \ref{lem4}, every unital positive linear map 
$\Phi: \cS\rightarrow M_n$ is completely positive.

\medskip

Next, we consider the case when $A$ is not normal.

Suppose $A \in M_2$. We may replace $A$ by
$U^*\left(e^{it}(A-(\tr A)I/2)\right)U$ for a suitable unitary $U\in M_2$,
and assume that 
$A = \begin{pmatrix} \alpha & b \cr 0 & -\alpha\cr\end{pmatrix}$
with $\alpha \ge 0$ and $ b>0$. Then $W(A)$ is an elliptical disk 
with major axis $[-r,r]$ and minor axis $i[-b,b]$
with $r = \sqrt{\alpha^2+b^2}$.
We may further apply an affine transform 
$$A = A_1  + i A_2 \mapsto \frac{1}{r} A_1 + \frac{i}{b}A_2.$$
Then  $A$ is unitarily similar to the 
symmetric matrix $C = \begin{pmatrix} i & 1 \cr 1 & -i\cr\end{pmatrix}$,
where $W(C)$ is the unit disk centered at the origin.
So, we may assume that
$\cS= {\rm span}\{I_2, C, C^*\} = {\rm span}\{E_{11}, E_{22}, E_{12}+E_{21}\}$
is the set of symmetric matrices in $M_2$.
By Lemma \ref{lem5},
every positive linear map $\Phi:\cS\to M_n$ is completely 
positive.

Finally,  suppose $A = A_0 \oplus [\mu]$. If $\mu\in W(A_0)$, 
then $W(A)=W(A_0)$ and the  result follows from the previous case. 
So we can assume that $\mu\not\in W(A_0)$.  
We may apply an affine transform to $A_0$ as in the preceding case
so that $W(A_0)$ is the unit disk centered at the origin,
and $A_0\in M_2$ is  nilpotent with norm 2.
Applying the same affine transform to $A$ 
yields $A = A_0 \oplus [\hat \mu]$.
Now, replacing $A$ by $e^{it}(U^*AU-\hat \mu I)$
for a suitable $t \in \IR$, we may assume that
$A = (rI_2 + C) \oplus [0]$, where $r = |\hat \mu| > 1$, and  
$C =  \begin{pmatrix} i & 1 \cr 1& -i \cr\end{pmatrix}$.
So, 
$A = \begin{pmatrix} r+i & 1 \cr 1 & r-i\cr\end{pmatrix} \oplus [0]$ 
with  $r > 1$.

\medskip
Suppose $B = H+iG  $ with $H = H^*$ and $G = G^*$.
We will construct  $P$ with $0 \le P \le I$ and a unital {\bf positive}
linear map
$\Psi: {\rm span}\{ E_{11}, E_{22}, E_{33}, E_{12}+E_{21}\}
\rightarrow M_n$ with
\begin{equation}\label{Psi}
\begin{array}{ll}
\Psi(E_{11}+E_{22}) = P, \quad
&\Psi(E_{12}+E_{21}) = H-rP, \\  &\\
\Psi(E_{11}-E_{22}) = G, \quad
&\Psi(E_{33}) = I-P.\end{array}
\end{equation}
By Lemma \ref{lem5}, $\Psi$ is completely positive. 
One easily checks that $\Phi$ is the restriction of $\Psi$
on $\cS$, and is also unital completely positive. 

\medskip
Let $t_0 \in (0, \pi)$ be such that $\cos(t_0) = -1/r$.
For any unit vector $\bx \in \IC^n$,
we have
\begin{eqnarray*}
\langle B\bx,\bx\rangle 
&\in&  W(B) \subseteq W(A) \\
 &=& \{a+ib: |b|\sqrt{r^2-1} \le a\mbox{ and } \\ 
  &&\hskip .7in 
     (a-r) \cos t + b \sin t \le 1
     \hbox{ for all } t \in [-t_0, t_0]\}.
  \end{eqnarray*}
Equivalently,
\begin{equation}\label{eq1a} \pm\sqrt{r^2-1} G\le H
\end{equation}
and
\begin{equation}\label{eq1b} \cos t(H-rI) +\sin t\, G\le I \quad
\hbox{ for all } t \in [-t_0,\ t_0].
\end{equation}

\medskip
Now, the desired map $\Psi$ satisfying 
(\ref{Psi})  is positive if and only if
$$\Psi\left(\begin{pmatrix} \cos^2\theta & \cos \theta \sin\theta \cr 
\cos \theta \sin\theta  & \sin^2\theta\cr\end{pmatrix} \oplus [0]\right)
\hskip 1in \ $$
$$\ \ 
=\frac{1}{2}\left(\cos^2\theta (P+G) +2 \cos \theta \sin \theta(H-rP) +
\sin^2 \theta (P-G)\right)$$
is positive semi-definite for all $\theta \in \IR$,
equivalently,
\begin{equation}\label{eq2.3}
P \ge \cos t (H-rP)  + \sin t \,G    \quad \hbox{ for all } t 
\in [-\pi,\pi].
\end{equation}
By (\ref{eq1a}),  $H\ge 0$ and  there exists a  
contraction $C=C^* \in M_n$ such that 
$G=\frac{1}{\sqrt{r^2-1}}H^{1/2}CH^{1/2}$. First, we show that
for $Q = \frac{I}{r+1} +\frac{C^2}{r^2-1}$,
\begin{equation}\label{eq3}
Q \ge
\cos t\left( I- rQ \right)+ \sin t \frac{C}{\sqrt{r^2-1}} =
\cos t\left( \frac{ I}{r+1} -\frac{rC^2}{r^2-1}\right) + \sin t \frac{C}
{\sqrt{r^2-1}}.
\end{equation}
To see this,  apply a unitary similarity to $C$ and 
assume that $C = \diag(c_1, \dots, c_n)$ with $c_j \in [-1,1]$. 
Then by Cauchy-Schwarz inquality and the fact that
$c_j^2 \in [0,1]$, we have 
\begin{eqnarray*}
&& \cos t\left( \frac{I}{r+1} -\frac{rc_j^2}{r^2-1}\right) 
+ \sin t \frac{c_j}{\sqrt{r^2-1}}  
\le  \sqrt {\left( \frac{1}{r+1} -\frac{rc_j^2}{r^2-1}\right)^2 
+\frac{c_j^2}{ r^2-1 }} 
\\ 
&\le& \sqrt{
\left( \frac{1}{r+1} -\frac{rc_j^2}{r^2-1}\right)^2 
+\frac{c_j^2}{ r^2-1 } + \frac{c_j^2(1-c_j^2)}{r^2-1}
} = \left(\frac{1}{r+1} +\frac{c_j^2}{r^2-1}\right)
\end{eqnarray*}
for each $j = 1, \dots, n$.
Hence
(\ref{eq3}) holds, and for $K=\frac{H}{r+1} +\frac{H^{1/2}C^2H^{1/2}}{r^2-1}\ge 0$, we have
\begin{equation}\label{eq2.5} 
\cos t H + \sin t  G \le (1+r\cos t) K
\quad \hbox{ for all }  t \in [-\pi,\ \pi]. \end{equation}
Suppose $V$ is unitary such that 
$K = V^* \diag(d_1, \dots, d_n)V$ with $d_1 \ge \dots \ge d_n \ge 0$.
Let 
$$P = \min\{K,I\} =  V^*\diag(p_1, \dots, p_n)V
\quad \hbox{ with  }  p_j = \min\{k_j,1\} \quad  \hbox{ for }
 j=1, \dots, n.$$
Then for $|t|\le t_0$, it follows from (\ref{eq1b}) and 
(\ref{eq2.5})  that
$$ \cos t H + \sin t  G \le (1+r\cos t)\min \{I,K\}
\le (1+r\cos t)P.  $$
For $\pi\ge |t|> t_0$, we have $1+r\cos t<0$. 
Together with (\ref{eq2.5}), 
we also have 

$$\cos t H + \sin t  G   
\le (1+r\cos t) K   
\le (1+r\cos t) P.$$
Thus, 
\begin{equation}\label{eq2.6}\cos t H + \sin t  G
\le (1+r\cos t)P \quad \hbox{ for all } t \in [-\pi, \pi]. \end{equation}
Hence, (\ref{eq2.3}) holds, and the result follows. \qed

\bigskip\noindent
{\bf Acknowledgment}

We would like to thank the referee for some helpful comments and the reference for Lemma \ref{lem2}. 
Li is an affiliate member of the Institute
for Quantum Computing, University of Waterloo, and is an
honorary professor of the
Shanghai University. His research was supported by USA
NSF grant DMS 1331021, Simons Foundation Grant 351047,
and NNSF of China Grant 11571220. Part of this research was done
while the two authors were visiting the Institute for Quantum Computing
at the University of Waterloo in the Fall of 2017. The
authors would like to express their thanks to the support of the 
Institute.

\end{document}